\documentclass[12pt,thmsa]{article}
\usepackage{amsfonts}
\usepackage{amssymb}
\newtheorem{teo}{Theorem}[section]

\newtheorem{obs2}[teo]{Remark}

\newtheorem{tea}{Theorem}[subsection]

\newtheorem{no2}[teo]{Note}

\newtheorem{no3}[tea]{Note}

\newcommand{\GL}{{\rm GL}}

\newcommand{\Q}{{\mathbb{Q}}}

\begin{document}
\title{{\bf On the level $p$ weight $2$ case of Serre's conjecture 
}}
\author{Luis Dieulefait
\\
Dept. d'Algebra i Geometria, Universitat de Barcelona;\\
Gran Via de les Corts Catalanes 585;
08007 - Barcelona; Spain.\\
e-mail: ldieulefait@ub.edu\\
 }
\date{\empty}

\maketitle

\vskip -20mm

\begin{abstract}
This brief note only contains a modest contribution: we just fix some inaccuracies in the proof of the prime level
 weight $2$ case of Serre's conjecture given in [K], for the case of trivial character. More precisely, the modularity lifting result needed at a crucial step is the one for the case of a deformation corresponding to a $p$-adic semistable (in the sense of Fontaine) Galois representation attached to a semistable abelian variety, but in [K] it is applied  a lemma only valid for potentially crystalline representations. The completion is easy: both the modularity lifting result applied in [K] (the one for the case of an abelian variety with potentially good reduction) and the one we need (the one for an abelian variety with bad semistable reduction) can be found in [D1]: they follow easily from a combination of modularity liftings results \`{a} la Wiles and other arguments.      

\end{abstract}
\section{The proof}

Acknowledgement: This note is written as an answer to the following question that S. Yazdani made us one month ago: Where is semistability used in the proof of the prime level weight $2$ case of Serre's conjecture given in [K]? We hope this note clarifies
 this point.\\  

We will only discuss the following case of Serre's conjecture: Prime level, weight $2$, trivial character. First recall the ingredients in the proof of this case of Serre's conjecture (cf. [K], except that we will change part (4)):\\
1) The truth of the level $1$ case of Serre's conjecture, arbitrary weight, proved in [K].\\
2) The result of existence of compatible families proved in [D2].\\
3) The result of existence of minimal lifts proved in [D3] and [KW].\\
4) An adequate modularity lifting result (a result that we will take from [D1], where it was used to prove modularity of certain abelian varieties).\\

The combination of these results gives elementary a proof of the case of Serre's conjecture we want to prove: the idea is just to
 follow the principle of ``switching the residual characteristic" (applied in [KW], [D3], [K] and [D4]) from one odd prime characteristic to another (in what follows $p$ is assumed to be an odd prime and $q$ is also assumed to be another odd prime, because the case of weight $2$ level $2$ is already proved in [D3] and [KW]). More precisely, given a mod $p$ representation $\rho$ of prime conductor (i.e., Serre's level) $q$, Serre's weight equal to $2$ and trivial character (an odd, irreducible, $2$-dimensional representation of the full Galois group of $\Q$), combining (2) and (3) we obtain a compatible family of Galois representations with prime conductor $q$, which corresponds to an abelian variety of $\GL_2$ type $A$ with good reduction away from $q$ and bad semistable reduction at $q$, such that for a prime dividing $p$ (in the field  $E$ of coefficients of this family) the corresponding $p$-adic representation is a lift of $\rho$. By definition, in order to show the modularity of $\rho$ it is enough to prove the modularity of this family, i.e., the modularity of $A$.\\
 Taking a prime dividing $q$ in $E$ we consider the corresponding $q$-adic representation in the family: it is of course enough to show that this member of the family is modular. If we consider the corresponding residual mod $q$ representation, since it is unramified outside $q$ it follows from (1) that it is either  modular or reducible.\\
 
  So, the essential argument needed in order to conclude the proof is a result already applied in [D1] (see [D1], last paragraph) which follows easily from modularity lifting results \`{a} la Wiles, and from the Mumford-Tate uniformization of semistable abelian varieties:\\
 
 {\bf Proposition:} \rm Let $q$ be an odd prime and consider an abelian variety $A$ of $\GL_2$ type defined over $\Q$ with  bad semistable reduction at $q$. Then, if the reduction mod $q$ of one of the $q$-adic representations attached to $A$ is known to be reducible or modular, $A$ is modular.\\
 
 Two remarks:\\
 1- A similar statement, but replacing ``bad semistable reduction" by ``good reduction" is also proved in [D1], using modularity lifting results \`{a} la Wiles and also a dihedral lemma (proved following ideas of Ribet) and a result we knew from Breuil concerning crystalline deformations in the case where the residual representation is in the case of level $1$ fundamental characters (this modularity lifting result was extended in [DM] and [D4] to the case of crystalline lifts of higher weight). This is the argument applied in [K] (these two results, both the dihedral lemma and the result of Breuil, Mezard and Savitt appear in section 5 of [K] and they are invoked in the proof of the prime level weight $2$ case of Serre's conjecture given in [K]), however it is not correct to apply this argument in the above situation (i.e., to prove the above proposition) because the result of Breuil, Mezard and Savitt applied in ([D1] and) [K] is a result for potentially crystalline deformations, not for semistable deformations.\\
  2- The proofs of the above proposition and of the similar statement in the case of good reduction given in [D1] are valid for any odd prime: in [D1] they are applied at the prime $3$ just because the residual modularity (or reducibility) for the abelian varieties considered in that article was known only at this prime, but the proofs work for any odd prime, no special properties of the prime $3$ are used (only that it is large enough: $3>2$). \\
  
  Conclusion: Thus, the above proposition concludes the proof of the modularity of $A$ and thus of the modularity of $\rho$:\\
  
  {\bf Theorem}: \rm The prime level weight $2$ (and trivial character) case of Serre's modularity conjecture is true. \\

  Remark: we have also deduced in particular the modularity of abelian varieties of $\GL_2$ type defined over $\Q$ with semistable reduction
   at one prime only (prime conductor).\\
   
   Final remark: For the case of prime level, weight $2$, and non-trivial character, the deformation is potentially crystalline and the proof given in [K] applies (we thank D. Savitt for this remark).

\section{Bibliography}
[D1] Dieulefait, L., {\it Modularity of abelian surfaces with Quaternionic Multiplication}, Math. Res. Letters {\bf 10} (2003), 145-150
\newline
[D2] Dieulefait, L., {\it Existence of compatible families and new cases of the Fontaine-Mazur conjecture},
J. Reine Angew. Math. {\bf 577} (2004) 147-151
\newline
[D3] Dieulefait, L., {\it The level $1$ weight $2$ case of Serre's conjecture}, preprint, (2004); available at: http://arxiv.org/math.NT/0412099 
\newline
[D4] Dieulefait, L., {\it Improvements on Dieulefait-Manoharmayum and applications}, preprint, (2005); available at: http://arxiv.org/math.NT/0508163 
\newline
[DM] Dieulefait, L.,  Manoharmayum, J., {\it Modularity of rigid Calabi-Yau threefolds over $\Q$}, in ``Calabi-Yau Varieties and Mirror Symmetry", Fields Institute Communications Series, American Mathematical Society {\bf 38} (2003) 159-166
\newline
[K] Khare, C., {\it On Serre's modularity conjecture for $2$-dimensional mod $p$ representations of $G_Q$ unramified outside $p$}, preprint, (2005); available at: http://arxiv.org/math.NT/0504080
\newline
[KW] Khare, C., Wintenberger, J-P., {\it On Serre's reciprocity conjecture for 2-dimensional mod p representations of the Galois group of Q}, preprint, (2004); available at: http://arxiv.org/math.NT/0412076 
\newline

\end{document}